\documentclass [a4paper]{IEEEtran}
\usepackage{amsmath}
\usepackage{amsthm}
\usepackage{amssymb, amsfonts}
\usepackage{graphics}
\usepackage{graphicx}
\usepackage{lipsum}
\usepackage[framed,numbered,autolinebreaks,useliterate,bw]{mcode}
\usepackage[utf8]{inputenc}
\usepackage[caption=false,font=footnotesize]{subfig}

\renewcommand{\u}[1]{\mathrm{#1}}

\newcommand{\dsort}{\mathrm{dsort}}
\newcommand{\diag}{\mathrm{diag}}
\newcommand{\diagmat}{\mathrm{diagmat}}

\title{Symmetric QR Algorithm with Permutations}
\author{Aravindh Krishnamoorthy\thanks{\hrule}\thanks{\noindent Aravindh Krishnamoorthy is currently with Ericsson Modem Nuremberg GmbH working in the area of Wireless Communication. Contact e-mail: aravindh.krishnamoorthy@ericsson.com, aravindh.k@ieee.org.}}

\begin{document}
\maketitle

\begin{abstract}
In this paper, we present the QR Algorithm with Permutations that shows an improved convergence rate compared to the classical QR algorithm. We determine a bound for performance based on best instantaneous convergence, and develop low complexity methods for computing the permutation matrices at every iteration. We use simulations to verify the improvement, and to compare the performance of proposed algorithms to the classical QR algorithm.
\end{abstract}

%%%%%%%%%%%%%%%%%%%%%%%%%%%%%%%%%%%%%%%%
\section{Introduction}
Let $A$ be an NxN symmetric invertible matrix, spectral theorem states that the matrix can be diagonalised as follows:
\begin{equation}
A = V \Lambda V^H
\end{equation}
where $\Lambda$ is a matrix with eigenvalues along the diagonals and $V = [v_1, v_2, \ldots, v_N]$ is a matrix with corresponding eigenvectors $v_1, v_2, \ldots, v_N$.

QR algorithm is an iterative algorithm based on QR decomposition to find the eigenvalues and eigenvectors simultaneously and is as follows:
\begin{align}
	A_0 &= A \label{eqn:qr1}\\
	A_k &= Q_k R_k \label{eqn:qr2}\\
	A_{k+1} &= R_k Q_k \nonumber\\
			&= Q_k^H A_k Q_k \label{eqn:qr3}
\end{align}
Further, $A_\infty = \Lambda$, $Q_\infty = I$, and $V = Q_0 Q_1 ... Q_\infty$.

For ease of analysis, we require in equation (\ref{eqn:qr2}) that upper-triangular matrix $R$ of QR factorization to contain positive values along the diagonal. This can be easily arranged by multiplying the row containing the negative value by `-1' and compensating the matrix $Q$. Such upper-triangular matrix $R$ can be shown to be the upper-triangular Cholesky factor of the matrix $A^2$.

For an introduction to the QR algorithm and a proof for convergence, and for modified algorithms: QR algorithm based on Hessenberg form and QR algorithm with shifts, refer \cite{book:golub} or \cite{paper:olver}.

\subsection{Cholesky Iteration}
A closely related algorithm to QR algorithm is Cholesky iterations based on Cholesky decomposition, given as follows:

\begin{align}
	B_0 &= A^2 \label{eqn:ch1}\\
	B_k &= R_k^H R_k \label{eqn:ch2}\\
	B_{k+1} &= R_k R_k^H \label{eqn:ch3}
\end{align}

The sets $\{A_0, A_1, A_2, ...\}$ and $\{B_0, B_1, B_2, ...\}$ are equivalent sets, and $A_k$ is the principal square root of the matrix $B_k$. This relation can be used to study the convergence properties of QR algorithm through Cholesky iteration which is more tractable for analysis; for e.g. see \cite{paper:schatz}.

\subsection{Permutations}
A permutation matrix $P$ is a matrix containing a single `1' in every row and column, and zeroes elsewhere. Let $\mathbb{P}$ be the set of such matrices of order N, then the cardinality of set $\mathbb{P}$ is $|\mathbb{P}| = N!$

A symmetric permutation $PAP^H$ on matrix $A$ defined above, maintains the symmetry of the matrix. It can be easily verified that the permuted matrix has the same eigenvalues as the original matrix, and the eigenvectors are $PV$.
 
\subsection{Terminology}
The following special functions are used in this paper:

\begin{itemize}
\item $\dsort: R^N \rightarrow R^N$ takes a real vector of order N as input, and returns a vector with the elements of input sorted in descending order.

\item $\diag: C^{N\u{x}N} \rightarrow R^N$ takes a symmetric matrix of order N as input, and returns a vector containing the diagonal elements of the matrix in the order along the main diagonal.

\item $\diagmat: R^N \rightarrow R^{N\u{x}N}$ takes a real vector of order N as input, and returns a real matrix containing the elements of the input vector along the main diagonal, and zeroes elsewhere.
\end{itemize}

\vspace{2mm}

In the following section, we present the QR algorithm with permutations along with a theoretical bound on performance. In section \ref{sec:pm} we discuss computation of the permutation matrices, in section \ref{sec:ssetup} we describe the simulation setup used to verify the performance of the algorithms. In section \ref{sec:results} we discuss the performance results for simulation runs involving positive-definite and symmetric matrices, followed by concluding remarks in section \ref{sec:con}.

\section{QR Algorithm with Permutations}
\label{sec:qrp}

Applying a permutation $P$ on equations (\ref{eqn:qr1}, \ref{eqn:qr2}, \ref{eqn:qr3}) yields an iteration of QR algorithm as follows:
\begin{align}
	PA_k P^H &= Q_k R_k \label{eqn:qrpf1}\\
	A_{k+1} &= R_k Q_k \nonumber\\
			&= (P^H Q_k)^H A_k (P^H Q_k) \label{eqn:qrpf2}
\end{align}
showing that a permutation on the matrix $A$ has an additional factor $P^H$ in the applied orthonormal matrix.

The effect of applying a permutation matrix on convergence can be intuitively understood when the QR algorithm is viewed as power iteration algorithm with a complete set of orthonormal basis-vectors instead of a single vector. 

Let n-th column of the matrix $Q_k$ be $q_n = a_{n,1} v_1 + a_{n,2} v_2 + ... + a_{n,N} v_N$, the convergence rate of the n-th eigenvalue $\lambda_n$ depends upon $a_{n,n}/a_{n,n-1} \lambda_n/\lambda_{n-1}$. With some luck, permutation changes the values $a_{n,1}, a_{n,2}, ..., a_{n,N}$, such that it offers a better convergence compared to the original matrix. 

Assume that $P$ is a permutation matrix that can ensure an optimum arrangement of the basis-vectors of $Q$, then the permutation causes the maximum convergence in that QR algorithm's iteration. This optimal arrangement can also aid the convergence of repeated eigenvalues or eigenvalues close to each other. However, this convergence is only speeded up by a linear factor!

%%%%%%%%%%%%%%%%%%%%%%%%%%%%%%%%%%%%%%%%
Let $A$ be an NxN symmetric invertible matrix as defined above, the QR algorithm with permutations with a permutation matrix $P_k$ at each iteration is given as follows:
\begin{align}
	A_0 &= A \label{eqn:qrp1}\\
	P_k A_k P_k^H&= Q_k R_k \label{eqn:qrp2}\\
	A_{k+1} &= R_k Q_k \nonumber\\
			&= (P_k^H Q_k)^H A_k (P_k^H Q_k) \label{eqn:qrp3}
\end{align}
Further, $A_\infty = \Lambda$, $Q_\infty = I$, and
\begin{equation}
V = P_0^H Q_0 P_1^H Q_1 ... P_\infty^H Q_\infty
\end{equation}

In general, in the k-th iteration, the estimates of eigenvalue matrix $\Lambda_k$, eigenvector matrix $V_k$ can be given as follows:
\begin{align}
\Lambda_k &= \diagmat(\diag(A_k)) \\
V_k &= P_0^H Q_0 P_1^H Q_1 ... P_k^H Q_k
\end{align}

The error in eigenvalue estimation $E_k$ for any algorithm can be given as follows: 
\begin{align}
E_k &= ||\dsort(\diag(\Lambda_k)) - \dsort(\diag(\Lambda))||_2 \label{eqn:ek}
\end{align}

QR algorithm with permutations requires the following additional computation compared to the classical QR algorithm: 1) computation of the permutation matrix, 2) symmetric permutation of matrix $A_k$ before QR decomposition, 3) permutation of the matrix $Q_k$ before multiplication, if eigenvectors are desired.

\subsection{Performance Bound}
The performance of the algorithm depends on the selection of the permutation matrices, i.e.  $E_k = f(A, P_0, P_1, ..., P_k)$.

A simple bound for $E_k$ can be given by requiring the best instantaneous convergence at every iteration. Therefore, the selection criterion for selecting the permutation matrix in the k-th iteration becomes:
\begin{equation}
	P_k = \arg \min_{\mathbb{P}} E_k^2 \label{eqn:bic}
\end{equation} 
That is, the permutation matrix $P_k$ of the available permutation matrices of order N, which when used in the k-th iteration, causes the least error in eigenvalue estimation $E_k$.

%%%%%%%%%%%%%%%%%%%%%%%%%%%%%%%%%%%%%%%%
\section{Computing the Permutation Matrices}
\label{sec:pm}

Computing the permutation matrix $P_k$ is non-trivial. In this section we present two simple methods for computing the permutation matrices which compute them without iterating through the set $\mathbb{P}$. Other methods for computing $P_k$ may be established easily, and their choice is a trade-off between convergence rate and computational complexity.

\subsection{Diagonal Ordering}
In this method, we choose the permutation matrix $P_k$ which orders the diagonal elements matrix $A_k$ in descending order of magnitude (absolute value). The MATLAB compatible pseudo-code to determine this permutation matrix is as follows:

\begin{lstlisting}
function P = diord(A) 
N = size(A,1) ;
D = abs(real(diag(A))) ;
[S,I] = sort(D, 'descend') ;
P = zeros(N,N) ;
P((0:N-1)'*N+I) = 1 ;
P = P' ;
\end{lstlisting}

\subsection{Column Ordering}
In this method, we choose the permutation matrix $P_k$ which orders the diagonal elements matrix $B_k = A_k^2$ in descending order. Such ordering also orders the columns of the matrix $A_k$ in descending order of their norm value.

This permutation matrix can be computed with a computational cost of $N^2$ for computing the N norm values of the columns of $A_k$. The MATLAB compatible pseudo-code to determine this permutation matrix is \mcode{P = diord(A*A)}.

\vspace{3mm}
Both diagonal ordering and column ordering methods do not preserve the shape of the matrix in subsequent QR iterations. Therefore, computational complexity reduction due to conversion to Hessenberg form cannot be applied with these two methods.

%%%%%%%%%%%%%%%%%%%%%%%%%%%%%%%%%%%%%%%%
\begin{figure*}[!t]
\begin{center}
  \mbox{\subfloat[]{\label{subfig:a} \includegraphics[width=\columnwidth,keepaspectratio=true]{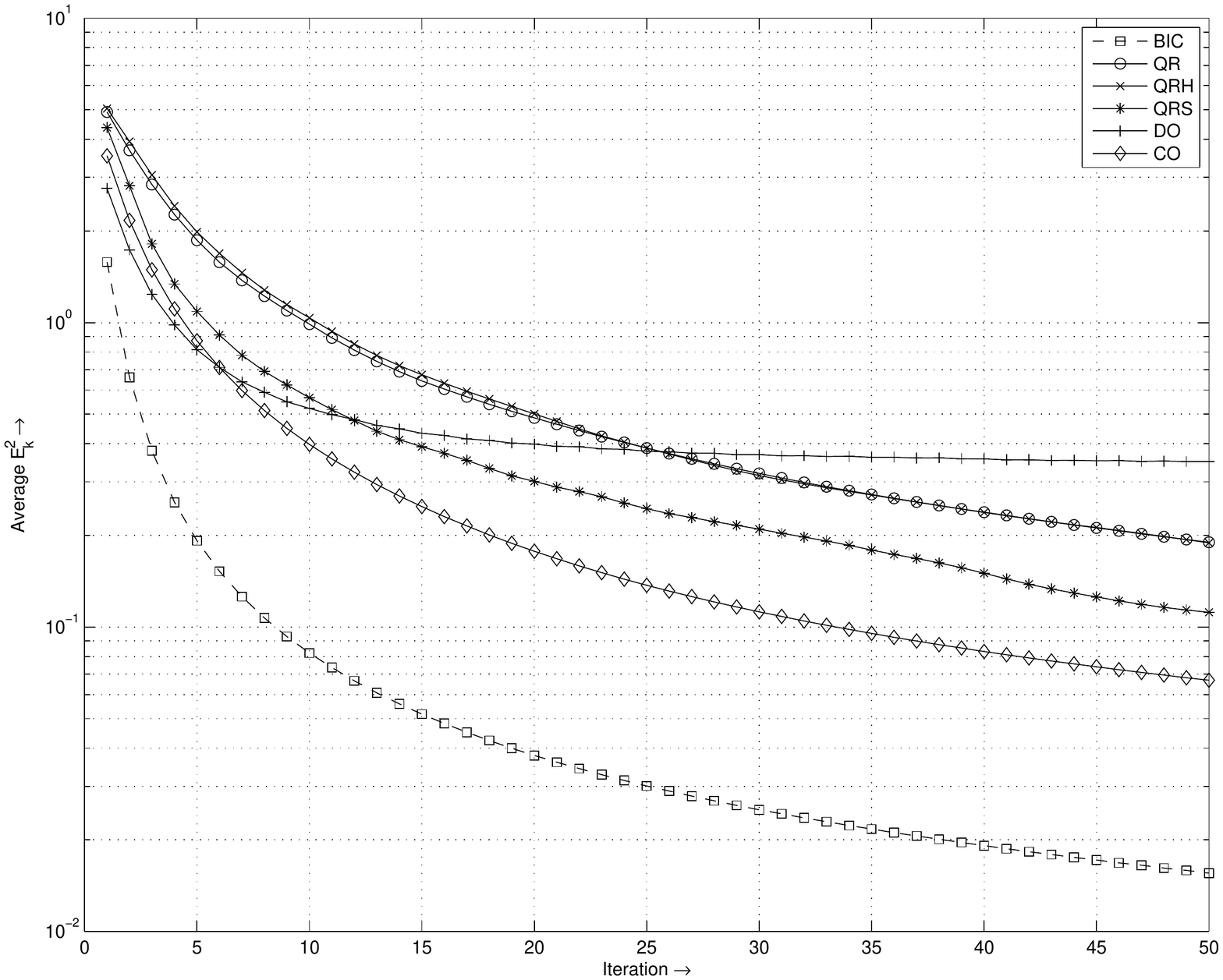}}}
  \mbox{\subfloat[]{\label{subfig:b} \includegraphics[width=\columnwidth,keepaspectratio=true]{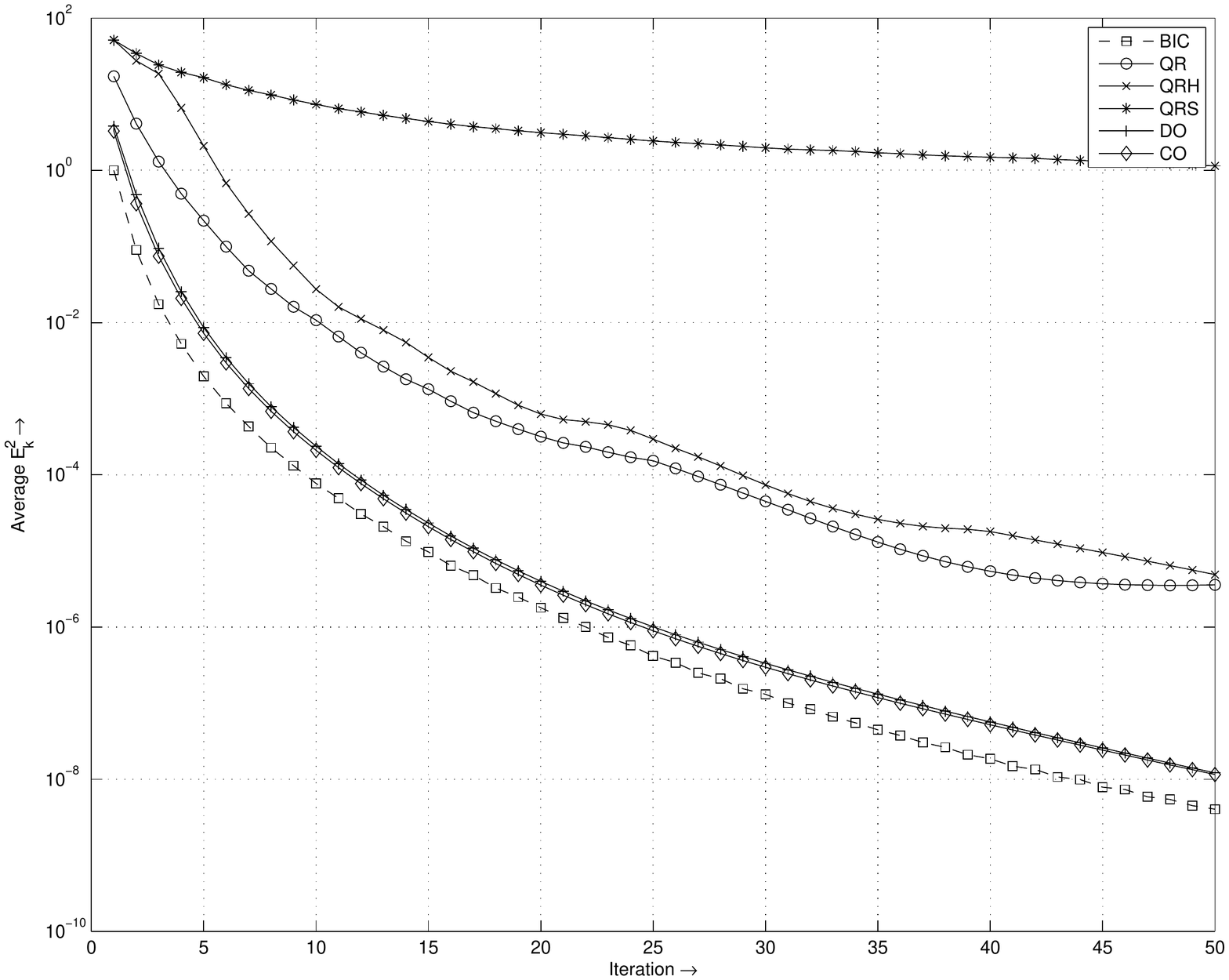}}}
\caption{Figures show the average eigenvalue estimation error $E_k^2$ for 4x4 symmetric invertible matrices, each with 50 iterations. (a) uses a random subset of size 10,000 of the set of symmetric matrices which includes positive-definite, negative-definite, and indefinite matrices. (b) uses a random subset of size 25,000 of the set of positive-definite symmetric matrices. The definition of legends is given in section \ref{sec:ssetup}.}
\label{fig:result}
\end{center}
\end{figure*}

\section{Simulation Setup}
\label{sec:ssetup}
In order to compare the performance of various algorithms, we use MATLAB based simulations. MATLAB compatible source code of the simulation is available from \cite{mex:code}.

The average eigenvalue estimation error $E_k^2$ equation (\ref{eqn:ek}) for 4x4 symmetric invertible matrices, each with 50 iterations are plotted. The first figure (\ref{fig:result}a) uses a random subset of size 10,000 of the set of symmetric matrices which includes positive-definite, negative-definite, and indefinite matrices. The second figure  (\ref{fig:result}b) uses a random subset of size 25,000 of the set of positive-definite symmetric matrices.

Label `QR' is the error plot for the classical QR algorithm, `QRH' first converts the given matrix into Hessenberg (tri-diagonal) form and then performs the classical QR algorithm, `QRS' uses the QRH algorithm ($A$ is in Hessenberg form) with shifts as given below ($A_0 = A$):
\begin{align}
	s_k &= A_{k-1}(N,N) \\
	A_k &= Q_k R_k - s_k I_N\\
	A_{k+1} &= R_k Q_k + s_k I_N
\end{align}
Here the notation $A_{k-1}(N,N)$ refers to the (N,N)-th element of the matrix $A_{k-1}$.

Label `DO' refers to the QR algorithm with permutations which uses diagonal ordering method for computing the permutation matrix at each iteration. `CO' uses the column ordering method.

Label `BIC' refers to the best instantaneous convergence bound as described in equation (\ref{eqn:bic}).

%%%%%%%%%%%%%%%%%%%%%%%%%%%%%%%%%%%%%%%%
\section{Simulation Results}
\label{sec:results}
From figure (\ref{fig:result}a), we observe that permutations increase the convergence rate on an average for symmetric and positive definite matrices.

The best instantaneous convergence (BIC) bound shows the highest rate of convergence, showing that permutation has an effect of speeding up the convergence of QR algorithm. However, meeting this bound using the definition is not feasible as it requires prior knowledge of eigenvalues.

QR algorithm with permutations with diagonal ordering (DO) shows good convergence rate for positive-definite matrices, but saturates after the first few iterations on an average, for symmetric matrices in general. DO performs better than classical QR algorithm and variations for positive-definite matrices.

QR algorithm with permutations with column ordering (CO) shows good convergence rate for positive-definite matrices and in general for symmetric matrices. CO performs better than classical QR algorithm and variations for positive-definite matrices and in general for symmetric matrices.

Therefore, for positive-definite matrices, DO provides the best convergence rate, and for symmetric matrices in general, CO provides the best convergence at low complexity for calculation of permutation matrices, both providing a 2x convergence speedup.

%%%%%%%%%%%%%%%%%%%%%%%%%%%%%%%%%%%%%%%%

\section{Conclusion}
\label{sec:con}
In this paper, we presented the QR algorithm with permutations that shows a significantly improved convergence rate compared to classical QR algorithm and variations. We determined a bound for performance based on best instantaneous convergence, and developed two low complexity methods for computing the permutation matrices at every iteration. We used simulations to verify that the performance bound is significantly improved upon using permutations, and compared the performance of QR algorithm with permutations using the two methods for computing permutation matrices.

We find that a 2x convergence speedup is obtained by using QR algorithm with permutations. Further, we find that the diagonal ordering provides the best trade-off for positive-definite matrices, and column ordering provides best trade-off for symmetric matrices in general.

\end{document}